\documentclass[12pt,peerreview,a4paper,onecolumn]{article}
\usepackage{graphicx}
\usepackage{latexsym}
\usepackage{amsfonts}
\usepackage{amsmath}
\usepackage{amssymb}
\usepackage{amscd}
\usepackage{amsxtra}
\usepackage{morefloats}
\usepackage{psfrag}
\usepackage{pstricks}
\usepackage{subfig}

\usepackage{amsxtra}
\usepackage{cite}
\usepackage[latin1]{inputenc}
\usepackage{setspace}
\usepackage{stfloats}\fnbelowfloat % Commands to control the presentation of floats.
\usepackage{morefloats}

\usepackage{multicol}
\usepackage{blindtext}

\newtheorem{assumption}{Assumption}
\newtheorem{theorem}{Theorem}

\newcommand{\realno}{\mbox{$I  \! \! R$}}
\newcommand{\naturalno}{\mbox{$I  \! \! N$}}

\date{}
\begin{document}
\title{Extended Hybrid Model Reference Adaptive Control of Piecewise Affine Systems}

\author{Mario di Bernardo \thanks{M. di Bernardo, U. Montanaro and  S. Santini are with
the Department of Electrical Engineering and Information Technology,
University of Naples Federico II,  Italy. E-mail: {\tt
\{mario.dibernardo, umberto.montanaro,
stefania.santini\}@unina.it}.} , Umberto Montanaro \thanks{U.
Montanaro is with the Department of Industrial Engineering,
University of Naples Federico II,  Italy. E-mail: {\tt
{umberto.montanaro}@unina.it}.} $^{,*}$ , Romeo Ortega
\thanks{R. Ortega is with the Laboratoire de Signaux et Systemes (SUPELEC), Paris, France. E-mail: {\tt
Romeo.Ortega@lss.supelec.fr}.} ,
\\ and Stefania Santini $^*$% <-this % stops a space
}

\maketitle

\begin{abstract}
This note presents an extension to the adaptive control strategy
presented in \cite{diMoSaTAC:13} able to counter eventual instability due to disturbances at the input of an otherwise  $\mathcal{L}_2$ stable closed--loop system. These disturbances are due to the presence of affine terms in the plant
and reference model.
%The extended strategy is able to address the possible
%instability caused by the presence of bounded $\mathcal{L}_\infty$ disturbances onto the
%closed loop system due to the presence of affine terms in the plant
%and reference model.
The existence of a common Lyapunov function is used to prove global convergence of the error system, even in the presence of sliding solutions, as well as boundedness of all the adaptive gains.
\end{abstract}

\section{Introduction}
As notably highlighted in \cite{Tao14},
adaptive control of switched systems is still an open problem. Recently, a novel model reference adaptive strategy has been
presented in  \cite{diMoSaTAC:13} that allows the control of
multi-modal piecewise linear (PWL) plants.
 Specifically, a hybrid model
reference adaptive strategy was proposed able to make a PWL plant
track the states of an LTI or PWL reference model even if the plant
and reference model do not switch synchronously between different
configurations. While stability is guaranteed for PWL systems, for
affine systems the presence of a non-square integrable disturbance
term in the error equations and the possible occurrence of sliding solutions can render the proof of
stability inadequate.

We wish to emphasize that the problem  of large state excursions and
instabilities caused by constant input disturbances on the closed
loop system is a common problem of adaptive control systems seldom
highlighted in the literature (see for example \cite{AsKa:07},
\cite{Anderson05}, and \cite{Anderson:1986} Sec. 4.4.4 p. 173).
Indeed, adaptive systems can be represented as the negative feedback
interconnection of a passive system (defined by the estimator) and a
strictly positive real (SPR) transfer function. A simple application
of the passivity theorem establishes that the overall system is
${\cal L}_2$-stable. However, this property does not ensure that the
system will remain stable in the presence of external disturbances
which are not ${\cal L}_2$. A simple example of this scenario is the
following system\cite{Ortega-pc}:
\begin{eqnarray*} \dot x_1 & = & f_1(x_1,x_2)=-x_1x_2\cos x_2\\ \dot x_2 & = & f_2(x_1,x_2)=x_1^2\cos x_2 - x_2+ \pi,
\end{eqnarray*}
where the presence of the constant input $\pi$ can cause exponentially growing trajectories that can significantly disrupt performance.

The aim of this note is to present a~
modification of the control strategy presented in
\cite{diMoSaTAC:13} able to guarantee asymptotic stability of the closed loop
system even in the presence of sliding mode trajectories and bounded $\mathcal{L}_\infty$
perturbations due to the affine terms in the description
of the plant and/or reference model. The idea is to add an extra
switching action to the controller in  \cite{diMoSaTAC:13} able to
compensate the presence of such a disturbance.
The proof of stability is obtained by defining an appropriate common Lyapunov function and analyzing
its properties along the closed-loop system trajectories within each of the phase space regions where the plant and reference model
are characterized by different modes, and along their boundaries.  We show that, even in the presence of possible sliding mode trajectories, the
origin of the closed-loop error system is rendered asymptotically stable by the extended strategy presented in this paper.
A preliminary version of the algorithm suitable to control bimodal piecewise affine system can be found in \cite{SIAM:2010}, \cite{dibeMoSA:08}, while experimental validation results are reported in \cite{CEP2012}.
A possible extension to discrete-time piecewise-affine (PWA) plants can be found in \cite{IJRNC:2013}, \cite{CEP:2013}.

\section{Problem statement and definitions}\label{sec:problem_and_defintion}
%Assume that the state space $\realno^n$ is partitioned into $M$
%polyhedral cells, say $\Omega_i$ ($i=0, \dots,M-1$), defined by
%generic hyperplanes of the form:
%\begin{equation}\label{eq:hyperplane_plant}
%    \Gamma_s: H_k^Tx+h_k = 0, \quad k = 1, \dots, \mathcal{K},
%\end{equation}
%for some $H_k \in \realno^n$ and $h_k \in \realno$.
%
%Furthermore, the switching surface between generic domains $\Omega_i$ and $\Omega_j$ is defined as:
%\begin{equation}\label{eq:hyperplane_Hij}
%\Gamma _{ij} :H_{ij}^T x + h_{ij}  = 0,
%\end{equation}
%with  $H_{ij}\in\realno^n$ and $h_{ij}\in \realno$. Note that the intersection is supposed to be not empty, {\em i.e.} $\bar \Omega_i \cap \bar \Omega_j\neq \emptyset$, being $\bar \Omega_i$  the closure of the generic cell $\Omega_i$.
Assume that the state space $\realno^n$ is partitioned by some
smooth boundaries into $M$ domains, say
$\left\{\Omega_i\right\}_{i\in\mathcal{M}}$ with $\mathcal{M} = \left\{0,1,\ldots, M-1\right\}$ such that
$\bigcup_{i=0}^{M-1}\Omega_i=\realno^n$ and, given two generic
indexes $i_1$ and $i_2\in \mathcal{M}$ (with $i_1\neq i_2$),  it follows $\Omega_{i_1}\cap\Omega_{i_2} = \emptyset$.

Let the plant be described by an $n$-dimensional multi-modal PWA system whose dynamics are given by:
\begin{equation}\label{eq:plant_model_pwa}
\dot x = {A_i}x + Bu + B_i\quad \text{if}\quad x \in {\Omega _i},
\quad i \in \mathcal{M},
\end{equation}
\noindent where $x \in \realno^n$  is the state vector, $u \in \realno$ is
the scalar input, and the matrices $A_i$, $B$, $B_i$ ($i =
0,1,\ldots, M-1$) are assumed to be in control canonical form, { i.e.}

\begin{equation}\label{eq:plant_matrices}
A_i  = \left[ {\begin{array}{*{20}c}
   0 & 1 &  \cdots  & 0  \\
   0 & 0 &  \ddots  &  \vdots   \\
    \vdots  &  \vdots  & {} & 1  \\
   {a_{i}^{(1)} } & {a_{i}^{(2)} } &  \cdots  & {a_{i}^{(n)} }  \\
\end{array}} \right], \;\;
B  = \left[
         \begin{array}{c}
           0 \\
           0 \\
           \vdots \\
           b \\
         \end{array}
       \right], \;\;
       B_i  = \left[
         \begin{array}{c}
           0 \\
           0 \\
           \vdots \\
           b_i \\
         \end{array}
       \right],
\end{equation}
with $b>0$. Note that all entries on the last row of the plant
matrices $A_i$,  $B$  and $B_i$ are supposed to be unknown. (Notice
that many generic bimodal PWL continuous systems can be transformed
into such a form as shown in \cite{TAC:2011}.)

The problem is to find an adaptive piecewise feedback law $u(t)$ to
ensure that the state variables of the plant track asymptotically
the states, say $\widehat x(t)$, of a reference model independently from their initial conditions.

Here, we assume that the reference model can be either an  LTI
system, or a multi-modal PWA system:
\begin{equation}\label{eq:reference_model_pwa}
\dot {\widehat x} = {{\widehat A}_{\widehat i}}\widehat x + \widehat
Br+\widehat B_{\widehat i} \quad \text{if}\quad \widehat x \in
{{\widehat \Omega }_{\widehat i}},\quad \widehat i \in
\mathcal{\widehat M},
\end{equation}
%\noindent where the state $\widehat x\in\realno^n$, $r \in \realno$ is the input to the reference model
%and the state space is partitioned into $\widehat M$ polyhedral cells, say $\widehat
%\Omega_{\widehat i}$ ($\widehat i=0, \dots,\widehat M-1$), through generic hyperplanes as:
%\begin{equation}\label{eq:hyperplane_refenrece_model}
%    \widehat \Gamma_{\widehat k}: \widehat H_{\widehat k}^Tx+\widehat h_{\widehat k} = 0, \quad {\widehat k} = 1, \dots,  \mathcal{\widehat K},
%\end{equation}
%being $\widehat H_{\widehat k} \in \realno^n$ and $\widehat h_{\widehat k} \in \realno$.
%%$\left\{\widehat \Omega_{\widehat i}\right\}_{\widehat i \in
%%\mathcal{\widehat M}}$ is a partition of $\realno^n$ with $\mathcal{\widehat M}\triangleq
%%\left\{0,1,\ldots, \widehat M-1\right\}$.
%
%Analogously we denote the switching surface between
%the domain $\widehat \Omega_{\widehat i}$ and $\widehat
%\Omega_{\widehat j}$, whose intersection is not empty, as:
%\begin{equation}\label{eq:hyperplane_Hij_hat}
%\widehat\Gamma _{\widehat i\widehat j} :\widehat H_{\widehat
%i\widehat j}^T x + \widehat h_{\widehat i\widehat j} = 0,
%\end{equation}
%with  $\widehat H_{\widehat i\widehat j}\in\realno^n$ and $\widehat
%h_{\widehat i\widehat j}\in \realno$.
\noindent where the state $\widehat x\in\realno^n$,
$\mathcal{\widehat M}\triangleq \left\{0,1,\ldots, \widehat
M-1\right\}$, $\left\{\widehat \Omega_{\widehat i}\right\}_{\widehat
i \in \mathcal{\widehat M}}$ is a partition of $\realno^n$ into
$\widehat{M}$ domains obtained by some smooth boundaries and $r \in
\realno$ is the input to the reference model. Note that the
reference model may possess a number of modes different from the one
of the plant, $\widehat M \neq M$. Furthermore, we assume that the
reference model defined as in \eqref{eq:reference_model_pwa} is
chosen so as not to exhibit sliding solutions and that it is
well-posed given the initial condition $\widehat x(0)=\widehat x_0$.
%Note that the absence of sliding is not a restrictive condition as it may appear at first.  For example, it is certainly satisfied by the large class of systems with degree of smoothness (DoS) 2.
In many practical cases, the aim of the control action can be that of compensating the discontinuous nature of the plant. In these situations, the control design presented above offers a simple and viable solution for this to be achieved by simply choosing a smooth or smoother reference model. This often corresponds to the conventional choice of an asymptotically stable LTI reference model in the case of smooth systems.

As for the plant, the matrices of the reference model are chosen to
be in the companion form given by ($\widehat i = 0,1,\ldots, \widehat M-1$):

\begin{eqnarray}\label{eq:modelref_matrices}
\widehat A_{\widehat i}  = \left[ {\begin{array}{*{20}c}
   0 & 1 &  \cdots  & 0  \\
   0 & 0 &  \ddots  &  \vdots   \\
    \vdots  &  \vdots  & {} & 1  \\
   {\widehat a_{\widehat i}^{(1)}} & {\widehat a_{\widehat i}^{(2)}} &  \cdots  & {\widehat a_{\widehat i}^{(n)} }  \\
\end{array}} \right],
\widehat B = \left[
               \begin{array}{c}
                 0 \\
                 0 \\
                 \vdots \\
                 \widehat{b}  \\
               \end{array}
             \right],
\widehat B_{\widehat i} = \left[
               \begin{array}{c}
                 0 \\
                 0 \\
                 \vdots \\
                 \widehat{b_{\widehat i}}  \\
               \end{array}
             \right]
             \end{eqnarray}
with $\widehat{b} >0$.

In what follows, we use the standard notation in
\cite{Branicky:98} (also adopted in \cite{HPassivity}), for both the switching instants of the plant and
reference model. More precisely, the switching sequence of the plant
is given by:

\footnotesize
\begin{equation}\label{eq:switching_sequence_plant}
\Sigma  = \{{\begin{array}{*{20}{c}}
{{x_0},\left( {{i_0},{t_0}} \right),\left( {{i_1},{t_1}} \right),
\left( {{i_2},{t_2}} \right) \ldots \left( {{i_p},{t_p}} \right) \ldots } \; |\; {{i_p} \in \mathcal{M},p \in \naturalno}  \\
\end{array}} \},
\end{equation}
\normalsize
where $t_0=0$ is the initial time instant and $x_0$ is the initial
state. Note that, as in \cite{Branicky:98}, when  $t\in[t_p\; ;
t_{p+1})$, then $x(t)$ belongs to $\Omega_{i_p}$ by definition and,
thus, the $i_p$-th subsystem is active. Obviously, the switching
sequence $\Sigma$ may be finite or infinite. If there is a finite
number of switchings, say $p$, then we set $t_{p+1}=\infty$.

For any $j\in\mathcal{M}$, we denote the sequence of switching times when the $j$-th subsystem is
switched on as:
\begin{equation}\label{eq:switching_s_j_plant_on}
\Sigma /j = \left\{ {\begin{array}{*{20}{c}}
   {{t_{{j_1}}},{t_{{j_2}}}, \ldots {t_{{j_s}}}, \ldots } & | & {{i_{{j_s}}} = j \; {\rm and} \;  s \in \; \naturalno}  \\
\end{array}} \right\},
\end{equation}
and, thus, the endpoints of the time intervals when the $j$-th subsystem is active can be given as:
\begin{equation}\label{eq:switching_s_j_plant_off}
\left\{ {\begin{array}{*{20}{c}}
   {{t_{{j_1} + 1}},{t_{{j_2} + 1}}, \ldots ,{t_{{j_s} + 1}}, \ldots } & | & {{i_{{j_s}}} = j \;   {\rm and} \;  s \in \naturalno}  \\
\end{array}} \right\}.
\end{equation}

Analogously,  we define the switching sequence of the reference model as:

\footnotesize
\begin{equation}\label{eq:switching_sequence_rm}
\widehat \Sigma  = \left\{ {\begin{array}{*{20}{c}}
 {{{\widehat x}_0},\left( {{{\widehat i}_0},{{\widehat t}_0}} \right),\left( {{{\widehat i}_1},{{\widehat t}_1}} \right),\left( {{{\widehat i}_2},{{\widehat t}_2}} \right) \ldots
 \left( {{{\widehat i}_p},{{\widehat t}_p}} \right) \ldots } \;| \; {{{\widehat i}_p} \in \mathcal{\widehat M},p \in \naturalno}  \\
\end{array}} \right\},
\end{equation}
\normalsize
with $\widehat t_0=0$. Hence, when  $t\in[\widehat t_p\; ; \widehat t_{p+1})$ then
$\widehat x(t)\in \widehat\Omega_{\widehat i_p}$ by definition and the $\widehat i_p$-th subsystem is active.

For any $\widehat j\in\mathcal{\widehat M}$ the sequence
of switching times when the $\widehat j$-th subsystem of the
reference model is switched on can be analogously defined as:
\begin{equation}\label{eq:switching_s_hat_j_rm_on}
\widehat \Sigma /\widehat j = \left\{ {\begin{array}{*{20}{c}}
   {{{\widehat t}_{{{\widehat j}_1}}},{{\widehat t}_{{{\widehat j}_2}}},
   \ldots {{\widehat t}_{{{\widehat j}_s}}}, \ldots } & | & {{{\widehat i}_{{{\widehat j}_s}}} = \widehat j,s \in \naturalno}  \\
\end{array}} \right\},
\end{equation}
with the endpoints of the intervals where the $\widehat j$-th mode is active being:
\begin{equation}\label{eq:switching_s_hat_j_rm_off}
\left\{ {\begin{array}{*{20}{c}}
   {{{\widehat t}_{{{\widehat j}_1} + 1}},{{\widehat t}_{{{\widehat j}_2} + 1}}, \ldots ,{{\widehat t}_{{{\widehat j}_s} + 1}},
    \ldots } & | & {{{\widehat i}_{{{\widehat j}_s}}} = \widehat j,s \in \naturalno}  \\
\end{array}} \right\}.
\end{equation}

We define the ''switching signals" $\sigma: \realno^+
\mapsto \mathcal{M}$ and $\widehat\sigma : \realno^+ \mapsto
\widehat{\mathcal{M}}$ as:
 \begin{eqnarray}
 \sigma (t) &=& i \quad \text{if}\quad x\left( t \right) \in {\Omega _i}, \quad
 \widehat\sigma (t) = \widehat i \quad \text{if} \quad \widehat x\left( t \right) \in {{\widehat \Omega }_{\widehat i}}  \label{eq:sigma_functions_all}
 \end{eqnarray}
and the indicator functions $\sigma_i(t)$ and $\widehat\sigma_{\widehat i}(t)$, as:
% \begin{eqnarray}
\begin{equation}\label{eq:sigma_function-a}
 {\sigma _i}\left( t \right) = \left\{
\begin{array}{l}
 1\; \text{if}\; x\left( t \right) \in {\Omega _i}, \\
 0\; \text{elsewhere}, \\
 \end{array} \right.
 \end{equation}
 \begin{equation}\label{eq:sigma_function-b}
{{\widehat \sigma }_{\widehat i}}\left( t \right) = \left\{
\begin{array}{l}
 1\; \text{if} \; \widehat x\left( t \right) \in {{\widehat \Omega }_{\widehat i}}, \\
 0\; \text{elsewhere}, \\
 \end{array} \right.%\label{eq:sigma_functions}
%\end{eqnarray}
 \end{equation}
with $ i = 0,1,\ldots,  M-1$ and $\widehat i = 0,1,\ldots, \widehat M-1$.
%Through the rest of the paper, we will assume that $\sigma$ and $\widehat{\sigma}$ have a finite number of discontinuities on every bounded time interval and take constant values on every interval between two consecutive switching times.

Also, $e_n\in \realno^n$ is defined as the basis vector
\begin{equation}\label{eq:Be}
e_n = \displaystyle \left[\begin{array}{cccc}
                            0 &  \ldots & 0 & 1
                          \end{array}
 \right]^T.
\end{equation}
%Finally, in the rest of the paper, we denote with the superscript $[v]$ the $v$-th component of a vector quantity of interest.

\section{Control Strategy}\label{sec:control_strategy}
The control problem described in Section
\ref{sec:problem_and_defintion} can be solved by means of an
extended switched adaptive strategy as described in the rest of this
section. The proposed approach extends the work presented in
\cite{diMoSaTAC:13} by exploiting an additional adaptive switching
control gain to cope with the presence of the bounded piecewise
constant input acting on the closed-loop system when the plant
and/or reference model are PWA.

\begin{assumption}\label{commonlyapunov}
Assume there exists a matrix  $P=P^T>0$ such that
\begin{equation}\label{eq:LMI_quad_stab_again}
 P\widehat A_{\widehat i}  + \widehat A_{\widehat i}^{T} P < 0 \quad  {\widehat i} = 0,1,2,\dots, {\widehat M}-1.
 \end{equation}
\end{assumption}

Given the above assumption our main result can be stated as follows.
%In this case it is possible to construct a Common Lyapunov Function to prove the following theorem.

\begin{theorem}\label{th:piecewise_aff}
Consider a PWA plant of the form \eqref{eq:plant_model_pwa} and a PWA reference model of the
form \eqref{eq:reference_model_pwa}. If the dynamic matrices ${\widehat A}_{\widehat i}$ of the reference model verify Assumption \ref{commonlyapunov}, then
the piecewise smooth adaptive control law: %in \eqref{eq:control_law}
\begin{equation}\label{eq:control_law}
    u(t)=K_R(t)r(t) + K_{FB}(t)x (t)+ K_A(t),
\end{equation}
where
\begin{eqnarray}
% \nonumber to remove numbering (before each equatio
K_R\left(t \right) &=& \alpha \int_0^t {y_e \left( \tau
\right)r\left( \tau  \right)d\tau  + \beta } y_e \left( t
\right)r\left( t \right),\label{eq:Kr}\\
K_{FB}\left(t \right) &=& K_0\left(t \right) + K_{\Sigma} \left(t
\right)+ \widehat K_{\Sigma}\left(t \right), \label{eq:Kfb}\\
K_{A}\left(t \right) &=& K_{0A}\left(t \right) + K_{\Sigma A} \left(t \right)+ \widehat K_{\Sigma A}\left(t \right)
\end{eqnarray}
with
\begin{equation}\label{eq:ye_xe}
y_e\triangleq C_{e}x_e , \quad x_e\triangleq \widehat x - x, \quad C_{e} \triangleq e_n^TP,
%\quad     \;\; \text{when} \;  \widehat
%x\left(t\right) \in \widehat \Omega_{\widehat i}, \;  \; \widehat i = 0,1,\ldots, \widehat M -1,
\end{equation}
%where matrices $C_{e\widehat i}$ are defined according to Assumption \ref{hyp:passivity}.
%and, in the case when Lemma \ref{hyp:passivity2} holds, they are set as
%\begin{equation}\label{eq:Cej}
%C_{e\widehat i}= e_n^T P_{\widehat i}.
%\end{equation}
\begin{eqnarray}
% \nonumber to remove numbering (before each equatio
K_0(t)&=& \alpha \int_0^t {y_e \left( \tau \right)x^T \left( \tau
\right)d\tau  + \beta } y_e \left( t \right)x^T \left( t \right),\label{eq:K0}\\
K_{\Sigma}\left(t \right) &=& \sum\limits_{j = 1}^{M - 1} {K_j
\left( t \right)}, \quad \widehat K_{\widehat \Sigma } \left( t
\right) = \sum\limits_{\widehat j = 1}^{\widehat M - 1} {\widehat
K_{\widehat j} \left( t \right)},\label{eq:K_sigma_hat_K_sigma}\\
K_j \left( t \right) &=& \left\{ \begin{array}{l}
 \rho \int_{t_{j_s } }^t {y_e \left( \tau  \right)
 x^T \left( \tau  \right)d\tau }, \quad \text{if} \quad x \in \Omega _j , \\
 0\quad \text{elsewhere},
 \end{array} \right.\label{eq:Kj}\\
 \widehat K_{\widehat j} \left( t \right) &=& \left\{ \begin{array}{l}
 \rho \int_{\widehat t_{\widehat j_s } }^t {y_e \left( \tau
 \right)x^T \left( \tau  \right)d\tau ,} \quad \text{if} \quad x \in \widehat \Omega _{\widehat j} , \\
 0\quad \text{elsewhere}, \label{eq:hat_Kj}
 \end{array} \right.,\\
K_{0A}(t)&=& \rho \int_0^t {y_e \left( \tau \right) d\tau}, \label{eq:K0A}\\
K_{\Sigma A}\left(t \right) &=& \sum\limits_{j = 1}^{M - 1} {K_{Aj} }, \quad
\widehat K_{\Sigma A} \left( t \right) = \sum\limits_{j = 1}^{M - 1} {{\widehat K}_{A {\widehat j}}{\widehat \sigma}_{\widehat j}}\label{eq:KhsA}
\\
K_{Aj} \left( t \right) &=& \left\{ \begin{array}{l}
 \rho \int_{t_{j_s } }^t {y_e \left( \tau  \right) d\tau }, \quad \text{if} \quad x \in \Omega _j , \\
 0\quad \text{elsewhere},
 \end{array} \right.\label{eq:KAj}\\
 \widehat K_{A\widehat j} \left( t \right) &=& \left\{ \begin{array}{l}
 \rho \int_{\widehat t_{\widehat j_s } }^t {y_e \left( \tau
 \right)d\tau ,} \quad \text{if} \quad x \in \widehat \Omega _{\widehat j} , \\
 0\quad \text{elsewhere}, \label{eq:hat_KAj}
 \end{array} \right.
\end{eqnarray}
and $\alpha$, $\beta$ and $\rho$ being some positive scalar constants, guarantees that the state tracking error $x_e(t)$ between the plant states $x(t)$ in \eqref{eq:plant_model_pwa} and the reference trajectory $\widehat x(t)$ in \eqref{eq:reference_model_pwa} converges asymptotically to zero, i.e. $\lim_{t \to \infty} x_e(t)=0$.
\end{theorem}

\subsection*{Remarks}
\begin{itemize}
%\item Note that hypothesis (\ref{eq:sliding_free_assumption}), or equivalently \eqref{eq:sliding_free_assumption_II}, is easy to be accomplished for systems of practical interest, as, for example, those described in \cite{MadiBeKe:00}.
\item The adaptation law presented above consists of three gains $K_R$, $K_0$ and $K_{0A}$ that remain switched on whatever the modes which
the plant and reference model are evolving in, together with some gains
$K_j$, $\widehat K_{\widehat j}$, $K_{Aj}$ and $\widehat K_{A {\widehat j}}$ that are switched on only when
the trajectories of the plant or reference model enter certain
domains in phase space. Specifically, the switching gains $K_j$ and $K_{Aj}$ are
associated to changes of the mode of the plant, whereas the
switching gains $\widehat K_{\widehat j}$ and $\widehat K_{A {\widehat j}}$ are associated to those of
the reference model. Furthermore, the gains $K_R$ and $K_0$ have the
same structure of the gains  in the Minimal Control Synthesis (MCS) approach \cite{PassivityMCS}, an application of Landau's Model Reference PI Adaptive Control scheme \cite{landau}.
%requiring no plant identification or off-line linear controller synthesis because of the canonical form assumed for  the plant and reference model \cite{StBe:90_1}.
\item  In order to compensate the bounded disturbance acting as an input onto the  closed-loop error system, the extended strategy exploits the additional adaptive term $K_A(t)$ when compared to the previous version of the algorithm presented in \cite{diMoSaTAC:13}, which is a set of switching integral actions used to properly compensate the affine term in each region.
\item At the generic $t_{j_s}$-th commutation,
the adaptive gains $K_j$ and $K_{Aj}$ are initialized to the last value assumed by
that gain when the trajectory of the plant $x(t)$ last exited from
region $\Omega_j$ (or zero otherwise). Analogously, the adaptive
gains $\widehat K_{\widehat j}$ and $\widehat K_{A {\widehat j}}$ at the generic $\widehat t_{\widehat
j_s}$~-~th commutation, is initialized with the last value assumed
by that gain when the trajectory $\widehat x(t)$ left the cell
$\widehat \Omega_{\widehat j}$ (or zero otherwise). Hence, according
to the notation used for the switching instants, we have:
\begin{eqnarray}
K_j \left( {t_{j_s } } \right) &=& K_j \left( {t_{j_{s - 1}  + 1} }
\right),s \ge 2,
\label{eq:reset_law_Kj}\\
K_{Aj} \left( {t_{j_s } } \right) &=& K_{Aj} \left( {t_{j_{s - 1}  + 1} }
\right),s \ge 2,
\label{eq:reset_law_KAj}\\
\widehat K_{\widehat j} \left( {\widehat t_{\widehat j_s } } \right)
&=& \widehat K_{\widehat j} \left( {\widehat t_{\widehat j_{s - 1} +
1} } \right),s \ge 2, \label{eq:reset_law_Kj_hat}\\
\widehat K_{A\widehat j} \left( {\widehat t_{\widehat j_s } } \right)
&=& \widehat K_{A\widehat j} \left( {\widehat t_{\widehat j_{s - 1} +
1} } \right),s \ge 2. \label{eq:reset_law_KAj_hat}
\end{eqnarray}
{Note that at the first transition the adaptive gains are set to zero, i.e.  $K_j \left( {t_{j_1 } } \right) = 0$, $K_j \left( {t_{j_1 } } \right) = 0$, $K_j \left( {t_{j_1 } } \right) = 0$, $\widehat K_{\widehat j} \left( {\widehat
t_{\widehat j_1 } } \right) = 0$, $\widehat K_{\widehat j} \left( {\widehat
t_{\widehat j_1 } } \right) = 0$. Furthermore, the integral part of the adaptive gains $K_R$ and $K_0$ in \eqref{eq:Kr} and \eqref{eq:K0} are set to zero at  time zero.}
\item {Both control gains $K_R$ and $K_0$ in \eqref{eq:Kr} and \eqref{eq:K0} have integral and proportional terms. It is worth remarking that the use of integral plus proportional adaptation has a beneficial effect upon the convergence of the generalized state error vector in comparison to the use of integral adaptation, specially at the beginning of the adaptation process \cite{landau}.  PI adaptation has also been used in \cite{AsKa:07}.}
\end{itemize}

\section{Proof of stability}
\noindent We now give the proof of Theorem \ref{th:piecewise_aff}
which is based on constructing an appropriate Common Lyapunov
Function (CLF) for the closed-loop system. Note that, due to the
presence of discontinuities in the closed-loop system dynamics, the
error state dynamics $x_e(t)$ is evaluated in the sense of Filippov
\cite{Fill:88} and hence, to prove convergence, the Lyapunov
function is also analyzed during possible instances of sliding
motion.

As the reference model \eqref{eq:reference_model_pwa} does not admit
sliding solutions by construction, we consider for the switched
closed-loop error system the following time-varying domains:
\begin{equation}\label{dominichiusi}
\left\{\Omega_i^c(t)\right\}_{i\in\mathcal{M}} =\{ x_e \in\realno^n :  \widehat{x}(t)-x_e \in \Omega_i \}_{i\in\mathcal{M}} ;
\end{equation}
with $\mathcal{M} = \left\{0,1,\ldots, M-1\right\}$. Furthermore we define as $\partial \Omega_i^{c}(t)$ their boundaries.  (Note that when the error trajectory $x_e(t)$ evolves along $\partial \Omega_i^{c}(t)$, the plant dynamics exhibit sliding motion along the surface $\partial \Omega_i$.)

From definition (\ref{dominichiusi}), it is now possible to express the indicator function $\sigma_i(t)$ in (\ref{eq:sigma_function-a}) as a function of the state tracking error as:
\begin{eqnarray}
 {\sigma _i}\left( t \right) &=& \left\{
\begin{array}{l}
 1\; \text{if}\; x_e\left( t \right) \in {\Omega_i^c(t)}, \\
 0\; \text{elsewhere}.\\
 \end{array} \right.
\label{eq:newsigma_functions}
\end{eqnarray}

From \eqref{eq:Kj}-\eqref{eq:hat_Kj} and
\eqref{eq:KAj}-\eqref{eq:hat_KAj}, it trivially follows that, at any
given time instant, only one of the pairs of adaptive gains $(K_1
K_{A1}), (K_2 K_{A2}),\ldots, (K_{M-1} K_{AM -1})$ and one of the
pairs  $(\widehat K_1 \widehat K_{A1}), \widehat (\widehat K_2
\widehat K_{A2}),\ldots,$ $(\widehat K_{\widehat M-1} \widehat
K_{A\widehat M-1})$  can be different from zero. Hence, equations
\eqref{eq:K_sigma_hat_K_sigma} and \eqref{eq:KhsA} can be easily
rewritten as:
\begin{equation}\label{eq:K_sigma_hat_K_sigma_2}
\begin{array}{l}
K_\Sigma  \left( t \right) = \sum\limits_{j = 1}^{M - 1} {\sigma _j \left( t \right)K_j \left( t \right)}, \quad \widehat K_{\widehat
\Sigma } \left( t \right) = \sum\limits_{\widehat j = 1}^{\widehat M - 1} {\widehat \sigma _{\widehat j} \left( t \right)\widehat
K_{\widehat j} \left( t \right)},\\
K_{\Sigma A}\left(t \right) = \sum\limits_{j = 1}^{M - 1} \sigma_j\left( t \right){K_{Aj} }, \quad
\widehat K_{\Sigma A} \left( t \right) = \sum\limits_{j = 1}^{M - 1} {\widehat \sigma} _{\widehat j} \left( t \right)
{\widehat K}_{A {\widehat j}},
\end{array}
\end{equation}
with $\sigma_j$  and $\widehat \sigma_{\widehat j}$ defined as in \eqref{eq:newsigma_functions} and \eqref{eq:sigma_function-b}, respectively.

Given expressions \eqref{eq:K_sigma_hat_K_sigma_2}, the adaptation law of the control gains \eqref{eq:Kj}, \eqref{eq:hat_Kj}, \eqref{eq:KAj},
\eqref{eq:hat_KAj} can be rewritten in terms of the indicator functions as ($i=1,\dots, M-1$, $\widehat i=1,\dots, \widehat M-1$):
\begin{eqnarray}\label{eq:new_swithing_gains}
% \nonumber to remove numbering (before each equatio
\dot K_i &=& \rho y_e(t)x^T(t) \sigma_i (t), \label{eq:dotKi}\\
\dot {\widehat K}_{\widehat i} &=& \rho y_e(t)x^T(t) \widehat \sigma_{\widehat i} (t), \label{eq:dotwidehatKi}\\
\dot K_{Ai} &=& \rho y_e(t) \sigma_i (t),\label{eq:dotKAi}\\
\dot {\widehat K}_{A \widehat i} &=& \rho y_e(t) \sigma_{\widehat i} (t).\label{eq:dotwidehatKAi}
\end{eqnarray}
Now, from the plant and the reference dynamics given in
\eqref{eq:plant_model_pwa} and \eqref{eq:reference_model_pwa}, respectively, by using the definition of the control strategy in
\eqref{eq:control_law}, after some algebraic manipulations the dynamics of the state tracking error can be rewritten as follows:
%\begin{equation}\label{eq:new_error_dynamics}
%\begin{array}{l}
%{\dot x_e} = {\dot {\widehat x}} - \dot x =
%{{\widehat A}_{\widehat \sigma(t)}}{x_e} + {e_n} \psi_I w -  {e_n}b\beta y_e w^Tw + e_n \sum\limits_{i = 1}^{M - 1} {\sigma_i {\psi _i} x}  \\+ e_n \sum\limits_{\widehat i = 1}^{\widehat M - 1} {{\widehat \sigma}_{\widehat i}{ \widehat \psi}_{\widehat i}} x + e_n ({\widehat b}_{\widehat \sigma (t)} - {b}_{\sigma (t)}) - e_nK_Ab ,
%\end{array}
%\end{equation}
\begin{equation}\label{eq:new2_error_dynamics}
\begin{array}{l}
{\dot x_e} = {\dot {\widehat x}} - \dot x =
{{\widehat A}_{\widehat \sigma(t)}}{x_e} + {e_n} \psi_I w -  {e_n}b\beta y_e w^Tw +\\
+ e_n \sum\limits_{i = 1}^{M - 1} {\sigma_i {\psi _i} x}  + e_n \sum\limits_{\widehat i = 1}^{\widehat M - 1} {{\widehat \sigma}_{\widehat i}{ \widehat \psi}_{\widehat i}} x + \\
+ e_n  \left [ \psi_{A0}  + \sum\limits_{{\widehat i} =1}^{{\widehat M}-1} ({\widehat \psi}_{A\widehat i} {\widehat \sigma}_{\widehat i} (t)) + \sum\limits_{{i} =1}^{{M}-1} ({\psi}_{Ai} {\sigma}_{i} (t))  \right  ], \end{array}
\end{equation}
where
%\footnotesize
\begin{eqnarray}\label{eq:gainspsi}
&& w \triangleq {\left[ {\begin{array}{*{20}{c}} {{x^T}}&{r}
\end{array}} \right]^T},\\
&&\psi_I \triangleq
\left[ {\begin{array}{*{20}{c}} {e_n^T\left( \widehat A_0 -
A_0 \right) - bK_0^I } & \vdots &{\widehat b - b{K_R^I}}
\end{array}} \right], \label{psiI}\\
&&{\psi _i} \triangleq \left[ {e_n^T\Delta {A_i} - b{K_i}}
\right], \; \Delta {A_i} \triangleq {A_0} - {A_i}, \label{eq:gainspsiA} \\
&&{{\widehat \psi }_{\widehat i}} \triangleq \left[ {e_n^T\Delta
{{\widehat A}_{\widehat i}} - b{{\widehat K}_{\widehat i}}} \right],
\; \Delta {{\widehat A}_{\widehat i}} \triangleq  {\widehat
A_{\widehat i}}-{\widehat A_0}, \label{eq:gainspsiB}\\
&&\psi_{A0}\triangleq ({\widehat b}_0 - {b}_0) - K_{0A}, \\
&&{\widehat \psi}_{A\widehat i} \triangleq \delta {\widehat b}_{\widehat i} - b{\widehat K}_{A\widehat i}, \; \delta {\widehat b}_{\widehat i} \triangleq ({\widehat b}_{\widehat i}- {\widehat b}_0),\\
&&{\psi}_{Ai} \triangleq \delta {b}_{i} - b{K}_{Ai}, \; \delta b_i \triangleq (b_i - b_0 ), \label{eq:gainpsilast}
\end{eqnarray}
%\normalsize
with $K_0^I, K_R^I$ being the integral part of $K_0$ in \eqref{eq:K0} and $K_R$ in
\eqref{eq:Kr}, respectively, and $i = 1, \ldots M-1$, $\widehat i = 1, \ldots \widehat M-1$.
Now, by means of the definition of the adaptive gains given in \eqref{eq:Kr}, \eqref{eq:K0}, \eqref {eq:K0A}, \eqref{eq:dotKi}, \eqref{eq:dotwidehatKi}, \eqref{eq:dotKAi}, \eqref {eq:dotwidehatKAi}, the dynamics of $\psi_I$, $\psi_i$,
${\widehat \psi}_{\widehat i}$, $\psi_{A0}$, ${\psi}_{Ai}$ and ${\widehat \psi}_{A\widehat i}$ in  \eqref{psiI}--\eqref{eq:gainpsilast} can be written as:
%\begin{eqnarray}\label{eq:gainsdotpsi}
%\dot \psi_I^T&=&-b \alpha y_e w, \label{eq:gainsdotpsiA}\\
%\dot{\psi}_{ i}^T &=& - b \rho y_e x {\sigma}_{i}\label{eq:gainsdotpsiB}\\
%\dot{\widehat \psi}_{\widehat i}^T &=& - b \rho y_e x {\widehat \sigma}_{\widehat i},\label{eq:gainsdotpsiC}\\
%{\dot{\psi}}_{A0}&=& - \rho y_e b, \label{eq:gainsdotpsiA1}\\
%{\dot{\psi}}_{Ai}&=& - \rho y_e b  {\sigma}_{i},\label{eq:gainsdotpsiA2}\\
%{\dot{\widehat \psi}}_{A\widehat i}&=& - \rho y_e b {\widehat \sigma}_{\widehat i}.\label{eq:gainsdotpsiA3}
%\end{eqnarray}
\begin{subequations}\label{clgains}
\begin{equation}
\dot \psi_I^T=-b \alpha y_e w, \label{eq:gainsdotpsiA}
\end{equation}
\begin{equation}
\dot{\psi}_{ i}^T = - b \rho y_e x {\sigma}_{i},\label{eq:gainsdotpsiB}
\end{equation}
\begin{equation}
\dot{\widehat \psi}_{\widehat i}^T = - b \rho y_e x {\widehat \sigma}_{\widehat i},\label{eq:gainsdotpsiC}
\end{equation}
\begin{equation}
{\dot{\psi}}_{A0}= - \rho y_e b, \label{eq:gainsdotpsiA1}
\end{equation}
\begin{equation}
{\dot{\psi}}_{Ai}= - \rho y_e b  {\sigma}_{i},\label{eq:gainsdotpsiA2}
\end{equation}
\begin{equation}
{\dot{\widehat \psi}}_{A\widehat i}= - \rho y_e b {\widehat \sigma}_{\widehat i}.\label{eq:gainsdotpsiA3}
\end{equation}
\end{subequations}
%
%Defining $\delta {\widehat b}_{\widehat i} = ({\widehat b}_{\widehat i}- {\widehat b}_0)$ and $\delta b_i = (b_i - b_0)$ ($i=1,\dots, M-1$, ${\widehat i}=1, \dots, {\widehat M} -1$), we have:
%\begin{eqnarray}\label{eq:bb}
%&&{\widehat b}_{\widehat \sigma (t)} = {\widehat b}_0 + \sum_{{\widehat i} =1}^{{\widehat M}-1} \delta {\widehat b}_{\widehat i}{\widehat \sigma}_{\widehat i} (t),\label{eq:bb1}\\
%&&{b}_{\sigma (t)} = \sum_{{i} =1}^{{M}-1} \delta {b}_{i}{\sigma}_{i} (t) .\label{eq:bb2}
%\end{eqnarray}
%
%According to \eqref{eq:bb1}-\eqref{eq:bb2}, the error system can be recast as:
%\begin{equation}\label{eq:new2_error_dynamics}
%\begin{array}{l}
%{\dot x_e} = {\dot {\widehat x}} - \dot x =
%{{\widehat A}_{\widehat \sigma(t)}}{x_e} + {e_n} \psi_I w -  {e_n}b\beta y_e w^Tw +\\
%+ e_n \sum\limits_{i = 1}^{M - 1} {\sigma_i {\psi _i} x}  + e_n \sum\limits_{\widehat i = 1}^{\widehat M - 1} {{\widehat \sigma}_{\widehat i}{ \widehat \psi}_{\widehat i}} x + \\
%+ e_n  \left [ \psi_{A0}  + \sum\limits_{{\widehat i} =1}^{{\widehat M}-1} ({\widehat \psi}_{A\widehat i} {\widehat \sigma}_{\widehat i} (t)) + \sum\limits_{{i} =1}^{{M}-1} ({\psi}_{Ai} {\sigma}_{i} (t))  \right  ] \end{array}
%\end{equation}
%where
%\begin{equation}\label{eq:psivariedef}
%\begin{array}{l}
%\psi_{A0}=({\widehat b}_0 - {b}_0) - K_{0A},\\
%{\widehat \psi}_{A\widehat i} = \delta {\widehat b}_{\widehat i} - b{\widehat K}_{A\widehat i},\\
%{\psi}_{Ai} = \delta {b}_{i} - b{K}_{Ai}.
%\end{array}
%\end{equation}
Note that, letting $z(t) \in
\realno^{\left(n+1\right)\left(\widehat{M}+M\right)-1}$ be the state
vector embedding the adaptive gain dynamics \eqref{clgains} as well
as the state tracking error \eqref{eq:new2_error_dynamics}
%$z(t):= [x_e \; \psi_I \; \psi_i \; {\widehat \psi}_{\widehat i} \; \psi_{A0} \;  {\psi}_{Ai} \;  {\widehat \psi}_{A\widehat i}]^T$,
the evolution of the closed-loop system \eqref{eq:new2_error_dynamics}, \eqref{clgains}  can be recast in a more compact form as the following set of differential equations with discontinuous right-hand side:
\begin{equation}\label{chiuso}
\dot z(t)=f_i(z) \quad x_e(t) \in \Omega_i^c \; i=0,\dots, M-1.
\end{equation}
where $f_i$ are the vector fields defined as the right-hand sides of the closed-loop system \eqref{eq:new2_error_dynamics},\eqref{clgains}.
\\
Now, let $P\widehat A_{\widehat i}  + \widehat A_{\widehat i}^{T} P =
-Q_{\widehat i}$ for some $Q_{\widehat i}=Q_{\widehat i}^T>0$,
${\widehat i} \in \mathcal{\widehat M}$ and let $V:
\realno^{\left(n+1\right)\left(\widehat{M}+M\right)-1}\to \realno$
be the Lipschitz, regular \cite{317122} and positive definite
candidate Lyapunov function given by:
\begin{equation}\label{lyapunov}
\begin{array}{l}
V=x_e^TPx_e + \frac{1}{\alpha b} \psi_I \psi_I^T + \frac{1}{\rho b} \sum\limits_{i = 1}^{M - 1} \psi_i \psi_i^T +\frac{1}{\rho b} \sum\limits_{\widehat i = 1}^{\widehat M - 1} {\widehat \psi}_{\widehat i}{\widehat \psi}_{\widehat i}^T +\\
+\frac{1}{\rho b} \psi_{A0}^2 + \frac{1}{\rho b} \sum\limits_{{i} =1}^{{M}-1} \psi_{Ai}^2 +\frac{1}{\rho b}\sum\limits_{{\widehat i} =1}^{{\widehat M}-1} {\widehat \psi}_{Ai}^2.
\end{array}
\end{equation}
To prove asymptotic stability of \eqref{chiuso} in what follows we will first evaluate of $\dot V$ in the interior of each generic region $\Omega_i^{c}$ and then along the generic surfaces $\partial S_l$ resulting from intersections of the manifolds $\partial \Omega_i^{c}$ (here $l=1,\dots, L$; $L$ being the number of manifolds where sliding is possible)  \cite{Alexander98}:
\begin{equation}\label{superficie}
\partial S_l = \bigcap_{d=1}^H \partial \Omega_{i_{d}}^c,
\end{equation}
with $H\leq (M-1)$.
%\\Note that classical Lyapunov stability theory cannot be directly applied as the right-hand side of \eqref{chiuso} is not Lipschitz \cite{RNC:RNC1317}.

\subsubsection{Evaluation of $\dot V$ in $\Omega_i^{c}$}
%Since the discontinuities may occur on switching surfaces $\partial S_l$, the Filippov set for each $x_e(t) \in \Omega_i^{c}$ involves only one elements and,
In the interior of each region $x_e(t) \in \Omega_i^{c}$, the error system \eqref{chiuso} is a smooth set of differential equations composed by equations \eqref{eq:new2_error_dynamics},\eqref{clgains} with $\sigma$ and $\widehat \sigma$ taking finite constant values associated to the active modes of the plant and reference model in that region.
The time derivative of $V$ along the trajectories \eqref{eq:new2_error_dynamics} can be computed as:

\footnotesize
\begin{equation}\label{lyapunovdott}
\begin{array}{l}
\dot V = -x_e^T Q_{\widehat \sigma(t)}x_e + \displaystyle \frac{2}{\alpha b} \psi_I {\dot \psi}_I^T + \frac{2}{\rho b} \sum\limits_{i = 1}^{M - 1} \psi_i {\dot \psi}_i^T + \frac{2}{\rho b} \sum\limits_{\widehat i = 1}^{\widehat M - 1} {\widehat \psi}_{\widehat i}{\dot {\widehat \psi}}_{\widehat i}^T +\\
+ 2 x_e^T P
 \left[  e_n \psi_I w - e_n b \beta y_e w^Tw
+ e_n  \sum\limits_{i = 1}^{M - 1} \sigma_i \psi_i x + \right.
 \left. e_n \sum\limits_{\widehat i = 1}^{\widehat M - 1} {\widehat \sigma}_{\widehat i} {\widehat \psi}_{\widehat i}x \right]+\\
 + 2x_e^TPe_n \left[ \psi_{A0} + \sum\limits_{{i} =1}^{{M}-1} ({\psi}_{Ai} {\sigma}_{i} (t)) + \sum\limits_{{\widehat i} =1}^{{\widehat M}-1} ({\widehat \psi}_{A\widehat i} {\widehat \sigma}_{\widehat i} (t)) \right]+\\
 \frac{2}{\rho b}{\psi}_{A0} {\dot \psi}_{A0} +\frac{2}{\rho b} \sum\limits_{{i} =1}^{{M}-1} {\psi}_{Ai} {\dot \psi}_{Ai}  +\frac{2}{\rho b} \sum\limits_{{\widehat i} =1}^{{\widehat M}-1} {\widehat \psi}_{A{\widehat i}} {\dot {\widehat \psi}}_{A{\widehat i}}.
\end{array}
\end{equation}

\normalsize
Substituting \eqref{eq:gainsdotpsiA}-\eqref{eq:gainsdotpsiA3} into equation \eqref{lyapunovdott} and taking into account that $x_e^TPe_n=e_n^TPx_e=y_e$, after some algebraic manipulations we have:
\begin{equation}\label{dotVdefneg}
\begin{array}{l}
\dot V=-x_e^TQ_{\widehat \sigma} x_e - 2 b \beta y_e^2 w^Tw  \leq  \\
\qquad\qquad \qquad - \min_{{\widehat i} \in \mathcal{\widehat M}}(\lambda_{\rm min}[Q_{\widehat i}])\|{x_e}\|^2 = - W(x_e),
\end{array}
\end{equation}
where ${\lambda_{\rm min}[Q_{\widehat i}]}$ is the smallest eigenvalue of $Q_{\widehat i}$.

\subsubsection{Evaluation of $\dot V$ along the manifolds $\partial S_l$}
In this case two different situations may occur; (i) the trajectory $x_e(t)$ crosses the generic manifold $\partial S_l$ \eqref{superficie} over a time interval of zero Lebesgue measure, or (ii) it exhibits sliding solutions. In the former case, the crossing has no effect on the stability analysis. Therefore, we focus below on the case where sliding occurs.\\
In particular, when sliding takes place, solutions should be interpreted in the sense of Filippov \cite{Fill:88}.
Using Filippov convex method, we consider the sliding vector field, say $f_F$, obtained by the convex combination \cite{DIECI}:
%Consider the generic $l$-th intersection manifold obtained intersecting some of the manifolds $\partial \Omega_i^{c}$ as \cite{Alexander98}:
%\begin{equation}\label{lyapunovdot}
%\partial S_l = \bigcap_{k=1}^H \partial \Omega_{i_{k}}^c,
%\end{equation}
%being $H\leq (M-1)$. Note that for $H=1$ we have a single switching surface, hence coinciding with one of the domain boundaries, $\partial \Omega_i^c$.
%The separating surface has co-dimension $H$, then (locally, in a neighborghood of $\partial S_l$) there are $2H$ regions and
%therefore $2H$ of the $M$ vector fields $f_i(z,t)$ in equation \eqref{chiuso}.
%Now, to characterize the sliding motion on switching surface $\partial S_l$, it is necessary to construct the following Fillippov (sliding) vector field lying in the tangent plane at $z \in \partial S_l$ \cite{DIECI}:
\begin{equation}\label{fF}
f_F:= \sum_{i=0}^{M-1}f_{i} (z) \gamma_{i}(z), \quad {\rm with} \;\gamma_{i}(z) \geq 0;
\end{equation}
where $f_i$ are the vector fields defined in equation \eqref{chiuso} and $\sum_{i=0}^{M-1} \gamma_{i}(z)=1$. \\Note that in the general case this is an underdetermined system of equations, hence there is no uniquely defined Filippov sliding vector. Since the following stability analysis does not depend on any particular choice of Filippov vector field, we do not consider this issue.\\
From \eqref{fF}, after some algebraic manipulations it is possible to write the closed-loop dynamics \eqref{eq:new2_error_dynamics},\eqref{clgains}
%\eqref{err-inc}
during the sliding motion as:
\begin{subequations}\label{sliding}
\begin{equation}
\begin{array}{l}
{\dot x_e}  =
{{\widehat A}_{\widehat \sigma(t)}}{x_e} + {e_n} \psi_I w -  {e_n}b\beta y_e w^Tw + e_n \sum\limits_{\widehat i = 1}^{\widehat M - 1} {{\widehat \sigma}_{\widehat i}{ \widehat \psi}_{\widehat i}} x \\
+  e_n \sum\limits_{{\widehat i} =1}^{{\widehat M}-1} {\widehat \sigma}_{\widehat i} {\widehat \psi}_{A\widehat i} + e_n  \psi_{A0}  +
e_n \sum\limits_{i= 1}^{M-1} \gamma_{i} ({\psi}_{i} x + {\psi}_{Ai})\label{err-incbn}
%PRR
%+ e_n \sum\limits_{i = 1}^{M - 1} {\sigma_i {\psi _i} x} \\ + e_n \sum\limits_{\widehat i = 1}^{\widehat M - 1} {{\widehat \sigma}_{\widehat i}{ \widehat \psi}_{\widehat i}} x +
%\left[ e_n  \left [ \psi_{A0}  + \sum\limits_{{\widehat i} =1}^{{\widehat M}-1} ({\widehat \psi}_{A\widehat i} {\widehat \sigma}_{\widehat i} (t)) + \sum\limits_{{i} =1}^{{M}-1} ({\psi}_{Ai} {\sigma}_{i} (t))  \right]  \right] ; \label{err-incb}
\end{array}
\end{equation}
\begin{equation}
\dot \psi_I^T=-b \alpha y_e w, \label{eq:gainsdotpsiAn}
\end{equation}
\begin{equation}
\dot{{\psi}_{ i}^T} = - \gamma_i b \rho y_e x,\label{eq:gainsdotpsiBn}
\end{equation}
\begin{equation}
\dot{{\widehat \psi}}_{\widehat i}^T = - b \rho y_e x {\widehat \sigma}_{\widehat i},\label{eq:gainsdotpsiCn}
\end{equation}
\begin{equation}
{\dot{\psi}}_{A0}= - \rho y_e b, \label{eq:gainsdotpsiA1n}
\end{equation}
\begin{equation}
{\dot{\psi}}_{Ai}= - \gamma_i b \rho y_e,\label{eq:gainsdotpsiA2n}
\end{equation}
\begin{equation}
{\dot{\widehat \psi}}_{A\widehat i}= - \rho y_e b {\widehat \sigma}_{\widehat i}.\label{eq:gainsdotpsiA3n}
\end{equation}
\end{subequations}
where only quantities depending on discontinuities due to plant switchings are convexified as the reference model is assumed not to exhibit sliding.

Evaluating the derivative of $V$ \eqref{lyapunov} (again in the sense of Filippov) along trajectories of \eqref{sliding}, some algebraic manipulations yield:
\begin{equation}\label{dotVdefneg2}
\dot V \leq - \min_{{\widehat i} \in \mathcal{\widehat M}}(\lambda_{\rm min}[Q_{\widehat i}])\|{x_e}\|^2 = - W(x_e).
\end{equation}

\subsubsection{Stability of the closed-loop adaptive system}
From \eqref{dotVdefneg}  and \eqref{dotVdefneg2}, it follows that for almost all $t$
\begin{equation}\label{barb1}
\dot V \leq - W(x_e) < 0 \; \; {\rm a.e.} \;t, \; \forall x_e \; \forall \psi_I, \psi_i, {\widehat \psi}_{\widehat i}, \psi_{A0}, {\psi}_{Ai}, {\widehat \psi}_{A\widehat i}.
\end{equation}
The derivative of the Lyapunov function along Filippov closed-loop solutions is negative, hence the origin of the closed-loop system is globally stable in the Filippov sense \cite{317122}.\\
Now following the approach in \cite{RNC:RNC1317}, from \eqref{barb1} for any closed-loop trajectory we have
\begin{equation}\label{barb2}
\sup_{t \in [ 0 + \infty)} V\leq C \quad \forall t \geq 0,
\end{equation}
with $C$ being a sufficient large positive constant. \\
From \eqref{barb1} and \eqref{barb2}, it follows
\begin{equation}\label{barb3}
\int^\infty_0 W(x_e(t))\,dt \leq C.
\end{equation}
Since $W(x_e)$ is a continuously differentiable positive-definite function, $W(x_e(t))$ is uniformly continuos.
Exploiting Barbalat's Lemma, $W(x_e(t))$ converges to $0$ as $t \to \infty$, hence the state tracking error $x_e(t)$ converges to $0$.
%PRRRRRRRRRRRRRRRRRRRR
%Hence $V$ is a common Lyapunov function \cite{Liberzon} of the
%switched closed-loop error system. Since $\dot V$ is definite
%negative then, from Lyapunov theory, the origin is a globally stable
%equilibrium point and then all the trajectories of the system and
%all adaptive gains are bounded (i.e. $x_e(t) \in \mathcal{L}_2 \cap
%\mathcal{L}_\infty$, $K_R$, $K_0$, $K_i$, ${\widehat K}_{\widehat
%i}$ $K_{0A}$, $K_{\Sigma A}$, ${\widehat K}_{\Sigma A}$ $\in \mathcal{L}_\infty$, independently from the initial
%condition). From \eqref{eq:control_law} and
%\eqref{eq:new2_error_dynamics}, it follows that $u(t)$ and $\dot
%x_e(t)$ are also bounded (i.e. $u(t), \dot x_e$ $\in
%\mathcal{L}_\infty$). Finally, as the evolution of all variables in
%the closed-loop system is bounded, Barbalat's Lemma can be used to guarantee
%asymptotic convergence of the state tracking error, $\lim_{t
%\to \infty} x_e(t)=0$ (see for example \cite{SaTa:10}).
$\Box$
\section{Conclusions}
We have presented an extension of the  hybrid adaptive strategy introduced in \cite{diMoSaTAC:13} aimed at compensating possible instabilities due to the presence of sliding mode trajectories and bounded perturbations acting on the closed loop error system. The external disturbances are due to the presence of affine terms in the plant and reference model dynamics. Using an appropriate common Lyapunov function, we have shown that the extended strategy guarantees asymptotic convergence of the tracking error, even in the presence of sliding solutions, as well as boundedness of all the adaptive gains.
%The pressing open problem remain of solving the disadvantage typical of MRAC approaches that,  as mentioned for example
%in \cite{Anderson:1986} Sec. 4.4.4 p. 173, in the presence of a
%bounded disturbance a drift of the adaptive gains in a model
%reference adaptive scheme can lead to ``exponentially growing
%transients that can significantly disrupt performance -- even if
%saturation is not reached''.

\bibliographystyle{unsrt}
\bibliography{MyBiblo_new}

\end{document}